\documentclass[11pt]{amsart}
\usepackage{graphicx}
\usepackage[centertags]{amsmath}
\usepackage{amsfonts}
\usepackage{amssymb}
\usepackage{amsthm}
\usepackage{newlfont}
\newtheorem{thm}{Theorem}[section]
\newtheorem{cor}[thm]{Corollary}
\newtheorem{lem}[thm]{Lemma}

\theoremstyle{definition}

\theoremstyle{remark}

\title{Almost conformal transformation in a class of Riemannian manifolds$^{1}$}
\author{Georgi Dzhelepov, Dimitar Razpopov, Iva Dokuzova}
\date{}
\begin{document}
\maketitle
\footnotetext[1]{This work is partially supported by project RS09 - FMI - 003 of the Scientific Research Fund, Paisii Hilendarski University of Plovdiv, Bulgaria}
%%%%%%%%%%%%%%%%%%%%%%%%%%%%%%%%%
\begin{abstract}

%%%%%%%%%%%%%%%%%%%%%%%%%%%%%%%%%%%
We consider a $3$-dimensional Riemannian manifold $V$ with a
metric $g$ and an affinor structure $q$.  The local coordinates of
these tensors are circulant matrices.  In $V$ we define an almost
conformal transformation. Using that definition we construct an
infinite series of circulant metrics which are successively almost
conformaly  related. In this case we get some properties.
\end{abstract}

\Small{\textbf{Mathematics Subject Classification (2010)}: 53C15,
53B20}

\Small{\textbf{Keywords}: Riemannian metric, affinor structure, almost conformal transformation
% ----------------------------------------------------------------

\section{Preliminaries}
\thispagestyle{empty}

We consider a $3$-dimensional Riemannian manifold $M$ with a
metric tensor $g$ and two affine tensors $q$ and $S$ such that:
their local coordinates form circulant matrices. So these matrices
are as follows:
\begin{equation}\label{f1}
    g_{ij}=\begin{pmatrix}
      A & B & B \\
      B & A & B \\
      B & B & A \\
    \end{pmatrix}, \quad A>B>0,
\end{equation}
where $A$ and $B$ are smooth functions of a point $p(x^{1}, x^{2},
x^{3})$ in some $F\subset R^{3}$,
\begin{equation}\label{f2}
    q_{i}^{.j}=\begin{pmatrix}
      0 & 1 & 0 \\
      0 & 0 & 1 \\
      1 & 0 & 0 \\
    \end{pmatrix},\qquad S_{i}^{.j}=\begin{pmatrix}
      -1 & 1 & 1 \\
      1 & -1 & 1 \\
      1 & 1 & -1 \\
    \end{pmatrix}.
\end{equation}

We note by $V$ the class of manifolds like $M$.

Let $M$ be in $V$ and $\nabla$ be the connection of $g$. Let us
give some results for $M$ in $V$, obtained in \cite{1}.
\begin{equation}\label{f3}
    q^{3}=E;\quad g(qu, qv)=g(u,v),\quad u,\ v\in \chi M .
\end{equation}
\begin{equation}\label{f4}
    \nabla q =0 \quad \Leftrightarrow \quad grad A=grad B . S .
\end{equation}
\begin{equation}\label{f5}
    0 < B < A \quad \Rightarrow \quad g \ is \ possitively \ defined.
\end{equation}
%%% ----------------------------------------------------------------

\section{Almost conformal transformation}

Let $M$ be in $V$.  We note
$f_{ij}=g_{ik}q_{j}^{k}+g_{jk}q_{i}^{k}$, i.e.
\begin{equation}\label{f2.1}
f_{ij}=\begin{pmatrix}
      2B & A+B & A+B \\
      A+B & 2B & A+B \\
      A+B & A+B & 2B \\
    \end{pmatrix}.
\end{equation}
We calculate $det f_{ij}=2(A-B)^{2}(A+2B)\neq 0$, so we accept
$f_{ij}$ for local coordinates of another metric $f$. Further, we
suppose $\alpha$ and $\beta$ are two smooth functions in $F\subset
R^{3}$ and we construct the metric $g_{1}$, as follows:
\begin{equation}\label{f2.2}
    g_{1}=\alpha .g +\beta .f.
\end{equation}
We say that equation (\ref{f2.2}) define an almost conformal
transformation, noting that if $\beta =0$ then (\ref{f2.2})
implies the case of the classical conformal transformation in $M$ \cite{2}.

From (\ref{f1}),(\ref{f2.1}) and (\ref{f2.2}) we get the local
coordinates of $g_{1}$:
\begin{equation}\label{f2.2*}
g_{1,ij}=\begin{pmatrix}
      \alpha A+2\beta B & \beta A+(\alpha+\beta)B & \beta A+(\alpha+\beta)B  \\
      \beta A+(\alpha+\beta)B  & \alpha A+2\beta B & \beta A+(\alpha+\beta)B \\
      \beta A+(\alpha+\beta)B  & \beta A+(\alpha+\beta)B & \alpha A+2\beta B \\
    \end{pmatrix}.
\end{equation}
We see that $f_{ij}$ and $g_{1,ij}$ are both circulant matrices.
\begin{thm}
Let $M$ be a manifold in $V$, also $g$ and $g_{1}$ be two metrics
of $M$, related by (\ref{f2.2}). Let $\nabla$ and $\dot{\nabla}$
be the corresponding connections of $g$ and $g_{1}$, and $\nabla
q=0$. Then $\dot{\nabla}q=0$ if and only if, when
\begin{equation}\label{f2.3}
    grad \alpha=grad \beta . S .
\end{equation}
\end{thm}
\begin{proof}
At first we suppose (\ref{f2.3}) is valid. Using (\ref{f2.3}) and
(\ref{f4}) we can verify that the following identity is true:
\begin{equation}\label{f2.4}
    grad (\alpha A+2\beta B)=grad (\beta A+(\alpha+\beta)B) . S
\end{equation}
The identity (\ref{f2.4}) is analogue to (\ref{f4}), and
consequently we conclude $\dot{\nabla}q=0$.

Inversely, if $\dot{\nabla}q=0$ then analogously to (\ref{f4}) we
have (\ref{f2.4}). Now (\ref{f4}) and (\ref{f2.4}) imply
(\ref{f2.3}). So the theorem is proved.
\end{proof}

Note. We see that (\ref{f2.4}) is a system of partial differential
equations. In this case we know that this system has a solution \cite{3}.

Let $w=w(x(p), y(p),z(p))$ be an arbitrary vector in $T_{p}M$,
$p\in M$, $M\subset V$, such that $qw\neq w$. For the metric $g$
of $M$ we suppose $0 < B < A$, i.e. $g$ is positively defined (see
(\ref{f5})).

Let $\varphi$ be the angle between $w$ and $qw$ with respect to
$g$. Then thank's to (\ref{f1}), (\ref{f2}) and (\ref{f3}) we get
    $cos\varphi=\dfrac{g(w, qw)}{g(w, w)}$, and we note that
    $\varphi\in (0, \dfrac{2\pi}{3})$ \cite{1}.
    \begin{lem}\label{1}
    Let $g_{1}$ be the metric given by (\ref{f2.2}). If $0 < \beta
    <\alpha$ and $g$ is positively defined, then $g_{1}$ is also
    positively defined.
    \end{lem}
    \begin{proof}
    For $g_{1}$ we have that $\alpha A+2\beta B- (\beta
    A+(\alpha+\beta) B= (\alpha-\beta)(A - B) > 0$. Analogously to
    (\ref{f2.1}) we state that $g_{1}$ is positively defined.
    \end{proof}
\begin{lem}\label{2}
    Let $w=w(x(p), y(p),z(p))$ be in $T_{p}M$,
$p\in M$, $M\subset V$, $qw\neq w$. Let $g$ and $g_{1}$ be the
metrics of $M$, related by (\ref{f2.2}). Then we have
\begin{align}\label{f11}
&g_{1}(w, w)=\alpha g(w, w)+2\beta g(w, qw)\\\nonumber
 &g_{1}(w,
qw)=\beta g(w, w)+(\alpha+\beta) g(w, qw).
\end{align}
    \end{lem}
\begin{proof}
    Using (\ref{f1}) and (\ref{f2}) we find
\begin{align}\label{f12}
&g(w, w)=A(x^{2}+y^{2} +z^{2})+2B(xy+yz+zx)\\\nonumber
 &g(w,
qw)=B(x^{2}+y^{2} +z^{2})+(A+B)(xy+yz+zx).
\end{align}
Now, we use (\ref{f2.2*}) and (\ref{f12}) after some computations
we get (\ref{f11}).
    \end{proof}
    \begin{thm}\label{4}
    Let $w=w(x(p), y(p),z(p))$ be a vector in $T_{p}M$,
$p\in M$, $M\subset V$, $qw\neq w$. Let $g$ and $g_{1}$ be two
positively defined metrics of $M$, related by (\ref{f2.2}). If
$\varphi$ and $\varphi_{1}$ are the angles between $w$ and $qw$,
with respect to $g$ and $g_{1}$ respectively, then the following equation is
true
\begin{equation}\label{f13}
\cos\varphi_{1}=\dfrac{\beta+(\alpha+\beta)cos
\varphi}{\alpha+2\beta cos\varphi}.
\end{equation}
    \end{thm}
\begin{proof}
    Since $g$ and $g_{1}$ are both positively defined metrics we
    can calculate $\cos \varphi$ and $\cos\varphi_{1}$, respectively \cite{2}. Then by
    using (\ref{f11}) from Lemma \ref{1} and Lemma \ref{2} we get
    (\ref{f13}).
    \end{proof}
    We note  $\varphi\in (0, \dfrac{2\pi}{3})$. Theorem \ref{4}
    implies immediately the assertions:
\begin{cor} If
$\varphi_{1}$ is the angle between $w$ and $qw$ with respect to
$g_{1}$ then $\varphi_{1}\in (0, \dfrac{2\pi}{3})$.
\end{cor}
\begin{cor}
Let $\varphi$ and $\varphi_{1}$ be the angles between $w$ and $qw$
with respect to $g$ and $g_{1}$. Then

1) $\varphi= \dfrac{\pi}{2}$ if and only if when $\varphi_{1}=
\arccos\dfrac{\beta}{\alpha}$ ;

2) $\varphi_{1}= \dfrac{\pi}{2}$ if and only if when $\varphi=
\arccos(-\dfrac{\beta}{\alpha+\beta})$.
\end{cor}
Further, we consider an infinite series of the metrics of $M$ in
$V$ as follows:
\begin{equation*}
    g_{0}, \ g_{1}, \ g_{2},\dots, \ g_{n}, \dots
\end{equation*}
where
\begin{equation}\label{f16}
   g_{0}=g,\quad g_{n}=\alpha g_{n-1}+\beta f_{n-1},\quad f_{n-1,is}= g_{n-1, ia}q_{s}^{a}+ g_{n-1,sa}q_{i} ^{a}, \quad 0 < \beta <
    \alpha.
\end{equation}
By the method of the mathematical induction we can see that the
matrix of every $g_{n}$ is circulant one and every $g_{n}$ is
positively defined.
\begin{thm}\label{5}
    Let $w=w(x(p), y(p),z(p))$ be in $T_{p}M$,
$p\in M$, $M\subset V$, $qw\neq w$. Let $\varphi_{n}$ be the angle
between $w$ and $qw$ with respect to metric $g_{n}$ from
(\ref{f16}). Then the infinite series:
\begin{equation*}
    \varphi_{0}, \ \varphi_{1}, \ \varphi_{2},\dots, \ \varphi_{n}, \dots
\end{equation*}
is converge and $\lim \varphi_{n}=0$.
    \end{thm}
\begin{proof}
    Using the method of the mathematical induction and Theorem \ref{4} we
    obtain
\begin{equation}\label{f18}
\cos\varphi_{n}=\dfrac{\beta+(\alpha+\beta)\cos
\varphi_{n-1}}{\alpha+2\beta \cos\varphi_{n-1}}
\end{equation}
as well as $\varphi_{n}\in (0, \frac{2\pi}{3})$. From (\ref{f18})
we get
\begin{equation}\label{f19}
\cos\varphi_{n}-\cos\varphi_{n-1}=\dfrac{\beta(1-\cos\varphi_{n-1})(1+2\cos
\varphi_{n-1})}{\alpha+2\beta \cos\varphi_{n-1}}.
\end{equation}
The equation (\ref{f19}) implies
$\cos\varphi_{n}>\cos\varphi_{n-1}$, so the series
$\{\cos\varphi_{n}\}$ is increasing one and since $\cos\varphi_{n}<
1$ then it is converge. From (\ref{f18}) we have
$\lim\cos\varphi_{n}=1$, so $\lim\varphi_{n}=0$.
    \end{proof}

\vspace{6mm}
\author{Georgi Dzhelepov\\ Department of Mathematics and Physics\\ Agricultural University of Plovdiv\\
12 Mendeleev Blvd.\\\vspace{6mm}Bulgaria 4000}\\
\author{Iva Dokuzova \\University of Plovdiv\\ FMI,  Department of
Geometry\\236 Bulgaria Blvd.\\Bulgaria 4003\\\vspace{6mm}
e-mail:dokuzova@uni-plovdiv.bg}\\
\author{Dimitar Razpopov \\ Department of Mathematics and Physics\\ Agricultural University of Plovdiv \\12 Mendeleev Blvd.\\
Bulgaria 4000 \\
e-mail:drazpopov@qustyle.bg}
% ------------------------------------------------------------------------------------------------------------------------
\end{document}